\documentclass{article}
\usepackage[utf8]{inputenc}
\usepackage[usenames]{color}
\usepackage{amsmath}
\usepackage{epsfig,latexsym}
\usepackage{epsf}
\usepackage{amsmath,amsfonts,amstext,amssymb,amsbsy,amsopn,amsthm,eucal}
\usepackage{latexsym}
\usepackage{epsfig,latexsym}
\usepackage[usenames]{color}
\usepackage{dirtytalk}
\usepackage{indentfirst}

\title{On the Philosophical Implications of the \\ Ouroboros Spaces and Their Functions}
\author{Nathan Thomas Provost$\footnote{Student of Applied Mathematics and Statistics at Brown University. University Email: nathan\_provost@brown.edu}$}
\date{}

\begin{document}

\maketitle

\begin{center}
\textbf{Abstract}
\end{center}

\small In this paper, I aim to articulate and investigate the philosophical implications and inherent symbolism surrounding the mathematical properties of Ouroboros spaces and their respective functions. Initially, I provide a brief historical background explaining how the symbol of the Ouroboros has been used and how it continues to be used as a term in mathematics. I then describe the philosophical symbolism and symbolic significance of the mathematical properties of the Ouroboros spaces and their functions, while offering an explanation as to why these concepts feel philosophically natural and intuitive. Following this discussion, I prove an aesthetically significant theorem that showcases the philosophical significance of the real Ouroboros functions. In closing, I articulate the interrelated, philosophical nature of these mathematical concepts, and describe how they impact other scientific fields both practically and philosophically.

\normalsize

\section*{The Philosophical and Historical Meanings of the Ouroboros in Mathematics}
\begin{center}
\say{\textit{There is repetition everywhere, and nothing is found only once in the world}.}
\end{center}
 --- Johann Wolfgang von Goethe\\

\normalsize Recently, I wrote a paper (\cite{1}) that proved two theorems regarding function spaces that consist of different kinds of \textit{Ouroboros} functions (functions of the form $f(f)=f$, as given in a paper I previously referenced by Soto-Andrade et al. \cite{2}). After investigating these \textit{Ouroboros} spaces, as I called them, I began to consider their place in the field of mathematical philosophy. The idea of a group of self-referencing, or alternatively, self-replicating objects seems natural in view of the world around us, which explains the intuitive concept behind the Ouroboros spaces and their functions. Yet, we must first ask ourselves: \textit{what is the Ouroboros and why is it important to the philosophy of mathematics}?\\

\newpage

Visually, the Ouroboros depicts a serpent or snake-like creature biting its own tail. Its origins lie in ancient Egypt, where it first became known as a symbol that represents the idea of a self-referential identity, symbolizing cyclical repetition and self-creation (as documented in \cite{3}). However, the symbol would persevere through the ages and leave a lasting impact on cultures around the world. A major scientific example occurred in the 19$^{th}$ century, when the German chemist A. K. von Stradonitz envisioned the Ouroboros in a dream. \cite{3} The ring-like shape of the serpent led him to determine the ring-like structure of the benzene molecule, which is still used today, but while this story is interesting, I still have not described the \textit{mathematical} significance of the Ouroboros.\\

The idea of a self-replicating equation has loosely existed in mathematics for many years (for example, the well-known expectation property $\mathbb{E}[\mathbb{E}[X]]=\mathbb{E}[X]$ seems to mirror the Ouroboros equation), but how did we decide to name $f(f)=f$, as given in \cite{2}, the Ouroboros equation? It first began when George Spencer-Brown, whose book, \textit{Laws of Form} \cite{4}, provided the foundation for a branch of logic/mathematics called the \textit{Calculus of Indications}. From here, Francisco Varela built upon Spencer-Brown's new subject in his paper titled \textit{A Calculus for Self-Reference} \cite{5}, the subject matter of which was of interest to his collaborative colleague Louis Kauffman. It is in Kauffman's work that the word \say{Ouroboros} is first used, as he describes Varela's development of Spencer-Brown's work as \say{the worm Ouroboros
embedded in a mathematical, non-numerical calculus.} \cite{6}. These ideas would later be discussed by Reichel in his paper \cite{7}, which quotes the same passage from Kauffman's work mentioning the Ouroboros. Much of these works discuss, or at least implicitly refer to cybernetics, which, in the broadest mathematical sense, focuses on systematic structures and structural information. In 2011, Soto-Andrade et al. labeled $f(f)=f$ the \say{Ouroboros equation} in their \say{mathematical exploration of self-reference and metabolic closure} \cite{2}, from which my previous work begins \cite{1}. These works each demonstrate the philosophical meaning of the Ouroboros by naming various forms of self-reference after the symbol. The idea of a structure that loops back to its own beginning is a very natural concept that stretches beyond mathematics alone (\textit{e}.\textit{g.} the water cycle, the life cycle, \textit{etc}.). Thus, it is only natural that we should explore the deeper philosophical meaning behind the mathematical concepts that carry on the Ouroboros' name.\\

While the works I have mentioned chiefly use the term Ouroboros to investigate cybernetics or logical structures, my goal was to provide an intuitive framework for the Ouroboros equation in functional analysis. For a univariate function $f(x)$, I called the function space:
\[\textbf{\textit{O}}(A)=\left\{f:A\rightarrow B \mid f(f(x))=f(x), \ \forall x\in A, \ \forall B\subseteq A \right\}\]
the Ouroboros space for the domain given by \textit{A} (which contains all the Ouroboros functions for the domain given by \textit{A}). \cite{1} More generally, I defined a higher dimensional Ouroboros space as follows:
\[\textbf{\textit{O}}(A^n)=\left\{f:A^n\rightarrow B \mid f(f(\textbf{x}),...,f(\textbf{x}))=f(\textbf{x}), \ \forall \textbf{x}\in A^n, \ \forall B\subseteq A\right\}\]
where $\textbf{x}=[x_1 \ ... \ x_n]^T\in A^n$ and $f(\textbf{x})=f(x_1,...,x_n)$. \cite{1} I plan to discuss the aesthetic value and philosophical symbolism behind this function space and its respective functions. Many of its properties reflect concepts that feel very natural to human beings, which embodies the spirit of the original symbol that carried great meaning in ancient Egypt. Furthermore, as we will see, Ouroboros functions embody unity in both the philosophical and mathematical sense, which further solidifies their aesthetic importance.

\section*{The Symbolism and Aesthetic Value of the \\Ouroboros Spaces and Their Functions}

Using the notation above, a univariate Ouroboros function must satisfy the equation $f(f(x))=f(x)$. This property truly embodies the ideas of repetition and self-replication. Furthermore, we can extend this property for countless iterations and the result is always the same:
\[f(f(f(...f(x)...)))=f(x)\]
This extension of the Ouroboros equation possesses potent philosophical meaning. For centuries, human beings have been preoccupied with the notion of preserving their legacy through their successors. In a way, the extended Ouroboros equation symbolizes the perpetuation of the self through generational replication. While this association is not entirely correct in the biological sense, it resonates with the fundamental concept of the cycle of life. The next successor born into a family acts as a replication of the self-same functions and their previous iterations. In a loosely philosophical manner, the concept of the replication of the self through continuous iterations extends to all animals.\\

More accurately, however, the aesthetic of this extended equation symbolizes the cycling of matter through the world around us. Water, for example, often appears on Earth in its liquid form, but moves from different states of matter. Yet, it returns to its former state in due time, no matter how many iterations of the cycle it goes through. Just as the original symbol depicts an endless loop of a snake biting its tail, the water cycle (and many other natural cycles for that matter) forever perpetuates itself in the form of an abstract, endless loop. Similarly, the concept of secondary succession, the environmental process through which a destroyed forest regrows, also echoes the symbolic sentiment of the Ouroboros equation. All of these cyclical concepts comprise a loose, metaphorical analog to the Ouroboros spaces, as they make up a collection of self-perpetuating entities that each exist under a set of formal conditions.\\

Yet, there is still a greater, unexplored aesthetic value that lies within the Ouroboros spaces. If we consider the univariate Ouroboros equation $f(f(x))=f(x)$ once more, and take the derivative of both sides of the equation, then we find, through the chain rule, that:
\[\frac{df(f(x))}{dx}=f'(f(x))\frac{df}{dx}=\frac{df}{dx} \ \therefore \ f'(f(x))=1\]
This result may not seem incredible, but the totality of the number 1 (which historically been referred to as \say{unity} \cite{8}) inspired me to investigate what I have come to call the \textit{Ouroboros Derivatives} (in this case, $f'(f(x))$ is the Ouroboros derivative in question). This exploration yielded the following aesthetic theorem.\\

\textbf{The Unified Ouroboros Derivatives Theorem}: Suppose that $f=f(x_1,...,x_n)=f(\textbf{x})$ is a continuous Ouroboros function of $n\geq1$ non-constant real variables where: 
\[f\in\textbf{\textit{O}}(\mathbb{R}^n)=\left\{f:\mathbb{R}^n\rightarrow B \mid f(f(\textbf{x}),...,f(\textbf{x}))=f(\textbf{x}), \ \forall \textbf{x}\in \mathbb{R}^n, \ \forall B\subseteq \mathbb{R}\right\}\] 
and $f:\mathbb{R}^n\rightarrow\mathbb{R}$. Let $f_{x_i}(\textbf{x})$ denote the partial derivative of $f$ with respect to $x_i$ for some $i\in\{1,...,n\}\subset\mathbb{N}$, which corresponds to the partial derivative of $f$ with respect to $g$ if $x_i=g$ (which we assume is also continuous and differentiable, along with the assumption that $f\in C^1(\mathbb{R}^n)$). Then, for any $n\in\mathbb{N}$, it holds that:
\[\sum_{i=1}^{n}f_{x_i}(f(\textbf{x}),...,f(\textbf{x}))=1 \ \ni \ f_{x_i}(f(\textbf{x}),...,f(\textbf{x}))=\frac{1}{n}, \ \forall i\in\{1,...,n\}. \]
where we refer to $f_{x_i}(f(\textbf{x}),...,f(\textbf{x}))$ as an Ouroboros derivative of $f(\textbf{x})$.\\

\noindent\textbf{Proof}: Suppose we have a function $f$ that meets the necessary specifications above. Let $x_i$ be an arbitrarily chosen variable of $f$. The multivariable Ouroboros equation states that $f(f(\textbf{x}),...,f(\textbf{x}))=f(\textbf{x})$. If we take the partial derivative with respect to $x_i$ of each side of the equation, we have:
\[\frac{\partial}{\partial x_i}f(f(\textbf{x}),...,f(\textbf{x}))=\frac{\partial f}{\partial x_i}\]
According to the multivariable chain rule, we have:
\[\frac{\partial}{\partial x_i}f(f(\textbf{x}),...,f(\textbf{x}))=\frac{\partial f}{\partial x_i}f_{x_1}(f(\textbf{x}),...,f(\textbf{x}))+ ... +\frac{\partial f}{\partial x_i}f_{x_n}(f(\textbf{x}),...,f(\textbf{x}))=\]
\[\frac{\partial f}{\partial x_i}\sum_{k=1}^{n}f_{x_k}(f(\textbf{x}),...,f(\textbf{x}))\]
Since $x_1=...=x_i=...=x_n=f(\textbf{x})$, it holds that $f_{x_1}(f(\textbf{x}),...,f(\textbf{x}))=...=f_{x_i}(f(\textbf{x}),...,f(\textbf{x}))=...=f_{x_n}(f(\textbf{x}),...,f(\textbf{x}))$, which means that:
\[\frac{\partial f}{\partial x_i}\sum_{k=1}^{n}f_{x_k}(f(\textbf{x}),...,f(\textbf{x}))=n\frac{\partial f}{\partial x_i}f_{x_i}(f(\textbf{x}),...,f(\textbf{x}))\]
Therefore, we now know that:
\[\frac{\partial}{\partial x_i}f(f(\textbf{x}),...,f(\textbf{x}))=n\frac{\partial f}{\partial x_i}f_{x_i}(f(\textbf{x}),...,f(\textbf{x}))\]
Returning to the derivative of the multivariable Ouroboros equation, we now have:
\[n\frac{\partial f}{\partial x_i}f_{x_i}(f(\textbf{x}),...,f(\textbf{x}))=\frac{\partial f}{\partial x_i}\]
which in turn means that:
\[f_{x_i}(f(\textbf{x}),...,f(\textbf{x}))=\frac{1}{n}, \ \forall i\in\{1,...,n\}\]
This subsequently holds for all variables $x_1,...,x_n$ since we selected $x_i$ arbitrarily. Naturally, it follows that:
\[\sum_{i=1}^{n}\frac{1}{n}=1\]
which ultimately means that:
\[\sum_{i=1}^{n}f_{x_i}(f(\textbf{x}),...,f(\textbf{x}))=1 \ \ni \ f_{x_i}(f(\textbf{x}),...,f(\textbf{x}))=\frac{1}{n}, \ \forall i\in\{1,...,n\} \ \qedsymbol\]
This theorem demonstrates the natural unity that exists within the Ouroboros spaces for real domains. Each of these Ouroboros derivatives are equal to one another, and beyond that, their sum is exactly one. The univariate example simply turned out to be a special case of this principle (when $n=1$). On a philosophical level, I must declare that all the Ouroboros derivatives of an Ouroboros function are \textit{intrinsically unified}. Just as the original symbol depicts the unification of a snake's head with its tail, the Ouroboros functions mirror this unity through their Ouroboros derivatives.\\

This idea of complete unity reflects the philosophical meaning behind many other mathematical and scientific  concepts. A continuous, Riemann-integrable probability density function must have an integral over all the real numbers that is equal to 1; $-e^{\pi i}=1$ (as a rearrangement of Euler's identity); and, on a fundamental level, 1 is the only number $n$ that satisfies the equation $nc=c$ for any scalar $c$. The last statement, the well-known multiplicative identity, reflects the philosophical idea of repetition and constancy, which appears beyond mathematics itself. The Law of Conservation of Mass states that matter cannot be created or destroyed in a closed system, and thus, the total amount of matter in such a system must remain the same. Indeed, the idea of \say{multiplying by 1} (equivalent to staying the same), or going through a process only to arrive at its initial state, is present throughout numerous branches of science and mathematics.
\newpage

\section*{Conclusion}

It is fascinating how a symbol from thousands of years ago can have such an influence on modern mathematics. From cybernetics, to functional analysis, the Ouroboros has inspired several mathematical concepts, whose aesthetic value and philosophical significance should not be ignored. The Ouroboros spaces offer a rigorous, mathematical description of a group of self-replicating entities, which is a concept that has existed philosophically for thousands of years. Our discussion of the symbolic nature of the Ouroboros equation demonstrated that, at its core, mathematics is a force that brings order to familiar concepts. The theorem I have proved herein further bolsters the aesthetic, philosophical power of the Ouroboros functions, as it ties their Ouroboros derivatives to the numerical representation of unity: 1. In further investigations, we should aim to investigate the relationship between Ouroboros spaces and probability theory, which I mentioned in my previous paper \cite{1}. It would also be interesting to explore the philosophical meanings and aesthetic value associated with this relationship. For now, however, I am satisfied with my philosophical exploration of the Ouroboros spaces and their functions, though I will undoubtedly pursue this concept further. Just as the snake biting its tail symbolizes a beginning emerging from an ending, the end of this paper marks the beginning of a new paper that has yet to be conceived.

\end{document}